\documentclass{article}

\PassOptionsToPackage{numbers}{natbib}

\usepackage{amsthm}
\usepackage{url}
\usepackage{algorithm}
\usepackage{algpseudocode}
\usepackage[breaklinks=true, hyperfootnotes=false]{hyperref}

\newtheorem{theorem}{Theorem}[section]

     \usepackage[final]{neurips_2021_ml4ps}

\usepackage[utf8]{inputenc} 
\usepackage[T1]{fontenc}    
\usepackage{hyperref}       
\usepackage{url}            
\usepackage{booktabs}       
\usepackage{amsfonts}       
\usepackage{nicefrac}       
\usepackage{microtype}      
\usepackage{xcolor}         

\usepackage{pgfplots}
\usepackage{tikz}
\usepackage[nottoc]{tocbibind}
\usepackage{float}
\usepackage{graphicx}
\usepackage{xcolor}
\usepackage{epsfig}
\usepackage{amssymb}
\usepackage{amsmath}
\usepackage{subfigure}
\usepackage{color}
\usetikzlibrary{arrows}
\usepackage{multirow}

\title{Modeling Advection on Directed Graphs using  Mat\'{e}rn Gaussian Processes for Traffic Flow}

\author{%
  Danielle C. Maddix
   \\
  Amazon Research\\
  2795 Augustine Dr. \\
  Santa Clara, CA 95054 \\
  \texttt{dmmaddix@amazon.com} \\
   \And
   Nadim Saad \thanks{Work conducted during an internship with Amazon AI.} \\
   Stanford University \\
   450 Serra Mall \\
   Stanford, CA 94305 \\
   \texttt{nsaad31@stanford.edu} \\
   \And 
   Yuyang Wang \\
   Amazon Research \\
   2795 Augustine Dr. \\
   Santa Clara, CA 95054 \\
   \texttt{yuyawang@amazon.com} \\
}

\begin{document}
\maketitle
\vspace{-0.5cm}
\begin{abstract}
The transport of traffic flow can be modeled by the advection equation.
Finite difference and finite volumes methods have been used to numerically solve this hyperbolic equation on a mesh.  
Advection has also been modeled discretely on directed graphs using the graph advection operator \cite{chapman2011, thesis}.  In this paper, we first show that we can reformulate this graph advection operator as a finite difference scheme. We then propose the Directed Graph Advection Mat\'{e}rn Gaussian Process (DGAMGP) model that incorporates the dynamics of this graph advection operator into the kernel of a trainable Mat\'{e}rn Gaussian Process to effectively model traffic flow and its uncertainty as an advective process on a directed graph.

\end{abstract}
\vspace{-.65cm}
\section{Introduction}
\vspace{-.25cm}
The continuous linear advection equation models the flow of a scalar concentration along a vector field. 
  The solutions to this hyperbolic partial differential equation may develop discontinuities or shocks over time depending on the initial condition. These shocks can model the formation of traffic jams, and their propagation along a road \citep{Richards1956}. 
Figure \ref{fig:upwind_exact} illustrates an example, where initially the first half of the road is 70\% occupied with cars, and the second half of the road is empty. The traffic propagates to the right until the whole road is 70\% occupied. 
Classical methods, such as finite differences and finite volumes, have been used to predict the flow of traffic along a road \cite{Lighthill1955OnKW, Richards1956}. 
These classical numerical methods do not incorporate any randomness into the model, and can be limited in incorporating the uncertainty among different driver's behaviors \cite{driving}.


\vspace{-.25cm}
 \begin{minipage}{\linewidth}
      \centering
      \begin{minipage}{0.49\linewidth} 
          \begin{figure}[H]
          \centering
              \includegraphics[width=\linewidth]{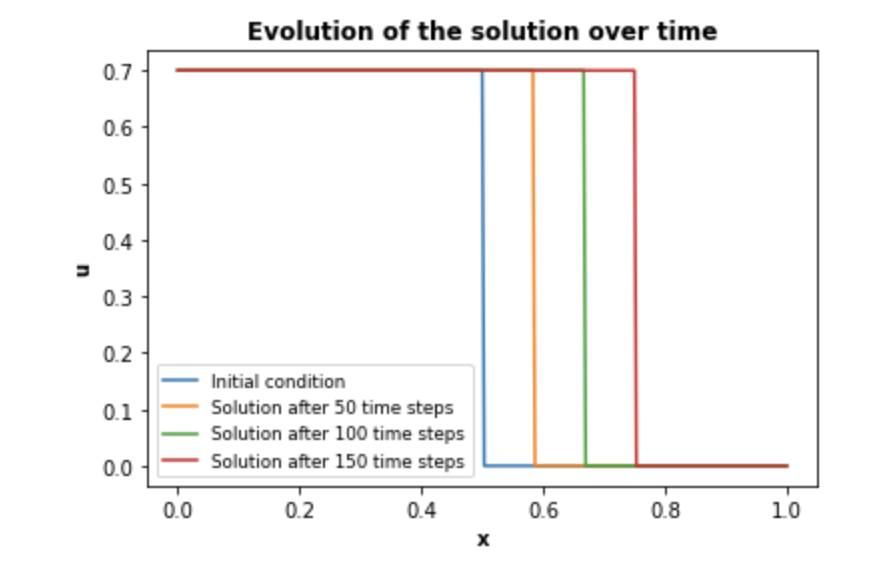}
              \caption{Propagation of cars on a road using an advection process.\\}
              \label{fig:upwind_exact}
          \end{figure}
      \end{minipage}
      \begin{minipage}{0.5\linewidth} 
      Gaussian processes (GPs) \cite{rasmussen2006} can learn unknown functions that allow use of prior information about their properties and for uncertainty modeling. \citet{kuper2020} propose the Gaussian Process Kalman Filter (GPKF) method to simulate spatiotemporal models, and test on the advection equation.  \citet{RAISSI2019686} train GPs on data to learn the underlying physics of non-linear advection-diffusion equations. Additional physics-based machine learning models \cite{borovitskiy2021matern} use the Matérn covariance function given below:
      \begin{equation} u \sim N\big(0,\big(\frac{2\nu}{\kappa^2} + \Delta \big)^{-\nu} \big), 
\label{eqn:Matern_GP}
\end{equation}
      \end{minipage}
  \end{minipage}
   where $u$ denotes an unknown function, $\nu < \infty$, $\kappa < \infty$ and $\Delta$ denotes the laplacian \cite{bakka2020diffusionbased}.  
 The Matérn kernel captures physical processes due to its finite differentiability, and is also commonly used to define distances between two points that are $d$ units distant from each other \cite{borovitskiy2021matern}.
 \citet{osti_1642956} propose training joint Mat\'{e}rn GPs to model space-fractional differential equations, in which the advection-diffusion equation is a special case.

Recent works including \cite{solomon2015pde} have studied solving partial differential equation (PDEs) on graphs. \citet{chapman2011, thesis} propose discrete advection and consensus operators to model advection and diffusion flows, respectively on directed graphs.  \citet{hosek2019} study the advection-diffusion equation on graphs using this discrete advection operator, and show that finite volume numerical discretizations can be reformulated as equations on graphs resulting in a corresponding maximum principle for this operator. Additional works have also looked at combining scientific computing and machine learning on graphs for spatiotemporal traffic modeling \cite{li2018}.  \citet{chamberlain2021grand} propose the Graph Neural Diffusion (GRAND) method, which combines traditional ODE solvers with graph neural networks (GNNs) to model diffusion on a undirected graph.  \citet{borovitskiy2021matern} propose to replace the continuous laplacian $\Delta$ in \eqref{eqn:Matern_GP} with the discrete graph laplacian operator $L$ to model diffusion on undirected graphs, which can be limited for traffic modeling. 

The goal of this paper is two-fold: to develop a model that effectively models traffic flow as an advective process on a directed graph and its uncertainty.  We propose a novel method, Directed Graph Advection Mat\'{e}rn Gaussian Process (DGAMGP) that uses a symmetric positive definite variant of the graph advection operator $L_{adv}$ as a covariance matrix in the Mat\'{e}rn Gaussian Process.  We use the square of the singular values of $L_{adv}$ to model the advection dynamics, and train a Mat\'{e}rn Gaussian Process to model the uncertainty.  We also show the connection between consistent finite difference stencils for solving the linear advection equation and the graph advection operator. Our novel linkage helps improve the understanding and interpretability of this graph advection operator. 

\section{Understanding the directed graph advection operator}
\label{sect:l_adv_prop}
We aim to model the continuous advection equation for unknown scalar $u$ under vector field $v$:
\begin{equation*}
        \frac{\partial u}{\partial t} = -\nabla \cdot (vu),
        \label{eqn:gov_eqn}
    \end{equation*}
stochastically on a directed graph. 
We define a directed, weighted graph $\mathcal{G} = (V, E, W)$ with $|V| = n$ nodes and $|E| = |W| = m$ edges, where $V$ denotes the vertex, $E$ the edge, and $W$ the edge weight 
sets, respectively. 
We discretize the flow $vu$ along edge $(i, j) \in E$ with weight $w_{ji} \in W$ as $w_{ji}u_i(t)$, where $u_i(t)$ denotes the concentration $u$ at node $i$ and time $t$. 

The graph advection operator $L_{adv}$ is defined so that the flow into a node equals the flow out of it \cite{chapman2011}:
      \begin{equation}
     \begin{aligned}
       \frac{du_i(t)}{dt} &= \sum_{j: (j, i) \in E} w_{ij}u_j(t) - \sum_{j: (i, j) \in E} w_{ji} u_i(t) = -[L_{adv}u(t)]_i,\\
    \end{aligned}       
    \label{eqn:graph_adv}
    \end{equation}
where $L_{adv} = D_{out} - A_{in}$ 
for diagonal out-degree matrix $D_{out}$ and 
in-degree adjacency matrix $A_{in}$.  
  For general directed graphs, $L_{adv}$ belongs to the square, non-symmetric with non-negative real part eigenvalues \cite{thesis} class of matrices in \cite{liesen2008}. 
   By design, $L_{adv}$ is conservative, unlike the related diffusion or consensus operator $L_{cons} = D_{in}-A_{in}$, where $D_{in}$ denotes the diagonal in-degree matrix \cite{chapman2011, thesis}.  
A main motivating reason for using $L_{adv}$ to model traffic flow is that it results in a conservative scheme. 

\vspace{-.25cm}
\paragraph{Reformulation of $L_{adv}$ as finite difference on balanced graphs.}
\label{subsec:finite_diff}
We notice that $L_{adv}$ at node $i$ is a weighted linear combination of the other nodes adjacent to it, which resembles finite difference stencils of the unknown and its neighbors.  We make this connection precise, and then construct example graphs where $L_{adv}$ corresponds to common finite difference schemes for linear advection.

\begin{theorem}
$L_{adv}$ corresponds to a semi-discrete finite difference advection scheme, where the sum of the coefficients is zero if and only if the graph $\mathcal{G}$ is balanced, i.e. $L_{adv} = L_{cons}$.
\end{theorem}
\vspace{-.5cm}
\begin{proof}
A finite difference approximation to the gradient can be written as the following weighted linear combination of its neighbors $u_j$ for arbitrary coefficients $c_{ij} \in \mathbb{R}$:
\begin{equation}
    -(u_x)_i \approx \sum_{j \ne i}c_{ij}u_j + c_{ii}u_i.
\label{eqn:FD}
\end{equation}
A consistent finite difference scheme is at least zero-th order accurate \cite{leveque}.   
Since the derivative of a constant is 0, the coefficients must sum to 0, i.e $c_{ii} = -\sum_{j \ne i}c_{ij}$. 
Combining 
\eqref{eqn:graph_adv} with \eqref{eqn:FD} gives:
$$(D_{out})_{ii} = \sum_{j: (i, j) \in E} w_{ji} = -c_{ii}  = \sum_{j \ne i}c_{ij} = \sum_{j: (j, i) \in E} w_{ij} = (D_{in})_{ii}.$$
The graph $\mathcal{G}$ is balanced by definition, and it follows that $L_{adv} = L_{cons}$. The other direction follows similarly.
\end{proof}

%

\begin{figure}[H]
\centering
\subfigure[first order upwind scheme]{\label{fig:graph_upwind}
\begin{tikzpicture}[->,>=stealth',shorten >=1pt,auto,node distance=2.5cm,
                    thick,main node/.style={circle,draw,font=\sffamily\bfseries}]
                    \centering

  \node[main node] (1) {$u_{i-1}$};
  \node[main node] (2) [right of=1] {$u_{i}$};
  \node[main node] (3) [right of=2] {$u_{i+1}$};

  \path[every node/.style={font=\sffamily\small}]
    (1) edge node[above] {$v/\Delta{x}$} (2)
    (2) edge node [above] {$v/\Delta{x}$} (3);
    \path (0,-1cm);
\end{tikzpicture}}
\subfigure[second order central scheme]{\label{fig:graph_central}
\begin{tikzpicture}[->,>=stealth',shorten >=1pt,auto,node distance=2.5cm,
                    thick,main node/.style={circle,draw,font=\sffamily\bfseries}]
\hspace{0.5cm}
  \node[main node] (1) {$u_{i-1}$};
  \node[main node] (2) [right of=1] {$u_{i}$};
  \node[main node] (3) [right of=2] {$u_{i+1}$};

  \path[every node/.style={font=\sffamily\small}]
    (1) edge node[above] {$v/2\Delta{x}$} (2)
    (3) edge node [above] {$-v/2\Delta{x}$} (2) 
    (2) edge [bend left] node[below] {$-v/2\Delta{x}$} (1)
    (2) edge [bend right] node[below] {$v/2\Delta{x}$} (3);
\end{tikzpicture}}
\caption{Balanced graphs on which $L_{adv}$ corresponds to finite difference stencils of linear advection.}
\end{figure}
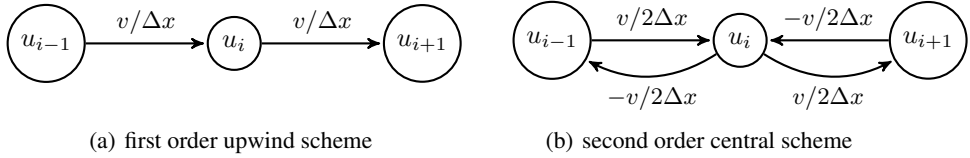

Applying $L_{adv}$ on the directed line graph in Figure \ref{fig:graph_upwind} results in the first order upwind scheme with spatial step size $\Delta x$ for $v > 0$ in \eqref{eqn:first_order_upwind} (See Appendix \ref{appendix:upwind} and Figure \ref{solve_upwind_conv} for the convergence study).
Similarly, Figure \ref{fig:graph_central} illustrates the directed graph in which $L_{adv}$ gives the second order central difference scheme, where $(u_x)_i \approx(u_{i+1}-u_{i-1})/(2 \Delta x)$ (See Appendix \ref{app:finite_diff_graphs} for additional examples).

\section{Directed Graph Advection Mat\'{e}rn Gaussian Process (DGAMGP)}
\label{sec:dgamgp}
We propose the novel Directed Graph Advection Mat\'{e}rn Gaussian Process (DGAMGP) model, which uses the dynamics of $L_{adv}$ to model advection stochastically on a directed graph through a discrete approximation to the continuous Laplacian $\Delta$ of the Mat\'{e}rn Gaussian Process in \eqref{eqn:Matern_GP}. The covariance matrix or kernel $\mathcal{K}$ of a Gaussian process needs to be symmetric and positive semi-definite.  This leads to some challenges with the $L_{adv}$ operator as it is not guaranteed in general to be symmetric or positive semi-definite (See Section \ref{sect:l_adv_prop}).  Note that using the graph Laplacian $L$ in the covariance matrix in the undirected graph case is more straightforward since $L$ is symmetric positive semi-definite. 

In our directed graph case,  we propose using $L^T_{adv}L_{adv}$ as the covariance matrix since it is symmetric positive definite, and hence orthogonally diagonalizable.  Analogous to \cite{borovitskiy2021matern}, we  define a function $\phi$ of a diagonalizable matrix through Taylor series expansion.   Then we can define its eigendecomposition as $L^T_{adv}L_{adv} = X_{adv}\Lambda_{adv} X^T_{adv}$, so that $\phi(L^T_{adv}L_{adv}) = X_{adv}\phi(\Lambda_{adv}) X^T_{adv}$, where $\phi(\Lambda_{adv})$ is computed by applying $\phi$ to the diagonal elements of $\Lambda_{adv}$.  

We compute the eigendecomposition of $L^T_{adv}L_{adv}
= V_{adv} \Sigma_{adv}^2 V_{adv}^T $, using the singular value decomposition (SVD) of $L_{adv} = U_{adv} \Sigma_{adv} V_{adv}^T$, 
where the eigenvalues and eigenvectors are the singular values squared and right singular vectors of $L_{adv}$, respectively. 
Hence, we model the advection dynamics using the square of the singular values of $L_{adv}$.
Our approach can also be viewed as adding the square of the singular values of $L_{adv}$ to the diagonal for regularization.  
Computing the thin-SVD is more computationally efficient and numerically stable, since we avoid explicitly forming the matrix-matrix product $L^T_{adv}L_{adv}$, which has double the condition number of $L_{adv}$, and the numerical issues with then computing its eigendecomposition.

We chose $\phi$ to be the Mat\'{e}rn covariance function in \eqref{eqn:Matern_GP}, and our DGAMGP model is given by:
 \begin{equation}
 u \sim N\big(0,\big( V_{adv} (\frac{2 \nu}{\kappa^2}I + \Sigma^2_{adv})^{-\nu} V^T_{adv}\big) \big).
  \label{eqn:directed_graph_gp}
 \end{equation}
This 
advective Gaussian Process is then trained on data by minimizing the negative log-likelihood of the Gaussian Process to learn the kernel hyperparameters $\nu$ and $\kappa$, and predict $u$ \cite{gardner2018}. For inference, we draw samples from the GP predictive posterior distribution with the learned hyperparameters \cite{rasmussen2006}.  See Algorithm \ref{alg:pseudocode} for details.

\paragraph{Choice of $L^T_{adv}L_{adv}$.} There are alternate approaches to symmetrize $L_{adv}$. The first simple approach explored is to utilize $L_{sym} = (L_{adv}^T + L_{adv}) / 2$ . This operator is not positive semi definite except in the balanced graph case. The second approach is to use the symmetrizer method in \cite{sen1988}, which generates a symmetric matrix $L'_{sym}$ with the same eigenvalues as $L_{adv}$ but is not always positive semi definite.
\begin{algorithm}[H]
\caption{The Directed Graph Advection Mat\'{e}rn Gaussian Process (DGAMGP)}\label{alg:pseudocode}
\begin{algorithmic}
    \item \textbf{Given} a directed graph $\mathcal{G} = (V,E, W)$ and training data $\mathcal{D} = \{(x_i, y_i)\}_{i=1}^n$.
    \begin{enumerate}

    \item Compute  $L_{adv}(\mathcal{G}) = D_{out} - A_{in}$.
    
    \item Compute the SVD of $L_{adv} = U_{adv} \Sigma_{adv} V_{adv}^T$.
    
    \item Generate a DGAMGP model in \eqref{eqn:directed_graph_gp}. 
    \item Minimize the GP negative log marginal likelihood using $\mathcal{D}$ to learn $\nu, \kappa$ and $\sigma$ \cite{gardner2018}.
    \item Given test data $\{x_i^*\}$, draw samples from the GP predictive posterior distribution \cite{rasmussen2006}.
\end{enumerate}
\end{algorithmic}
\end{algorithm}
\section{Numerical Results}
\label{sec:num_res}
\vspace{-0.25cm}
          In this section, we utilize our DGAMGP model for traffic modeling on synthetic and real-world directed traffic graphs.  
        The data $\mathcal{D} = \{(x_i, y_i)\}_{i=1}^n$ denotes the traffic flow speed in miles per hour $y_i$ at location $x_i$. We test our model's predictive ability to predict the velocities of cars on a road at different positions.  We use hold-out cross validation to split the data points generated into training (70\% of the data) and testing data (30\% of the data).  We extend the code in \cite{borovitskiy2021matern} to compute the singular value decomposition of $L_{adv}$ to train our DGAMGP model on a directed graph.  The code is available at
\url{https://github.com/advectionmatern/Modeling-Advection-on-Directed-Graphs-using-Mat-e-rn-Gaussian-Processes}, and the experiments are run on Amazon Sagemaker \cite{sagemaker}.

\begin{figure}[H]
\centering
\subfigure[intersection graph]{\label{fig:merging_lanes}
\begin{tikzpicture}[->,>=stealth',shorten >=1pt,auto,node distance=1.95cm,
                    thick,main node/.style={circle,draw,font=\sffamily\bfseries}]
  \node[main node] (1) {$u_{i-2}$};
  \node[main node] (2) [right of=1] {$u_{i-1}$};
  \node[main node] (3) [right of=2] {$u_i$};
  \node[main node] (4) [below of=1] {$u_{i-4}$};
    \node[main node] (5) [below of=2] {$u_{i-3}$};
    \node[main node] (6) [right of=3] {$u_{i+1}$};

  \path[every node/.style={font=\sffamily\small}]
    (1) edge node[above] {$v/\Delta{x}$} (2)
    (2) edge node [above] {$v/\Delta{x}$} (3) 
    (4) edge node[above] {$v/\Delta{x}$} (5)
    (5) edge [bend right] node[below] {\hspace{.5cm} $v/\Delta{x}$} (3)
    (3) edge node[above] {$2v/\Delta{x}$} (6);
\end{tikzpicture}}
\subfigure[loop graph]{\label{fig:loop}
\begin{tikzpicture}[->,>=stealth',shorten >=1pt,auto,node distance=2cm,
                    thick,main node/.style={circle,draw,font=\sffamily\bfseries}]
  \node[main node] (1) {$u_1$};
  \node[main node] (2) [right of=1] {$u_{i-1}$};
  \node[main node] (3) [right of=2] {$u_i$};
    \node[main node] (4) [right of=3] {$u_n$};

  \path[every node/.style={font=\sffamily\small}]
    (1) edge node[above] {$v/\Delta{x}$} (2)
    (2) edge node[above] {$v/\Delta{x}$} (3)
    (3) edge node[above] {$v/\Delta{x}$} (4)
    (4) edge [bend left] node[below] {$v/\Delta{x}$} (1);
\end{tikzpicture}}
\caption{Graphs representing two lanes merging into one (\textit{left}) and a loop (\textit{right}).}
\end{figure}
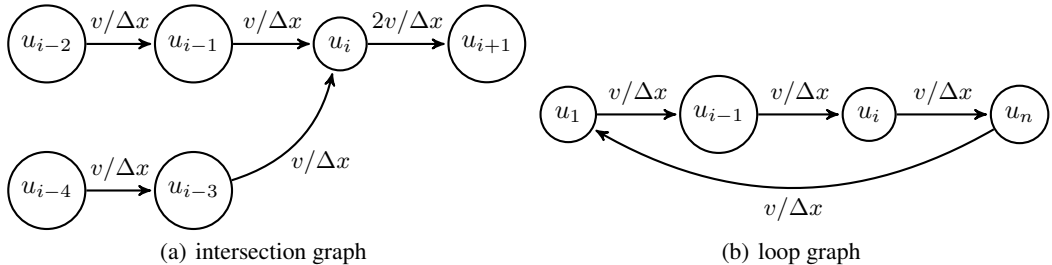

\paragraph{Regression results on synthetic graphs.}
 We generate synthetic data that models traffic along a road, which has a relatively high density of cars in the first half and a low density of cars in the second half. We train and test our model on the upwind scheme in Figure \ref{fig:graph_upwind}, central scheme in Figure \ref{fig:graph_central}, an intersecting lane graph, where two lanes merge into one lane in Figure \ref{fig:merging_lanes} and a loop graph representing the upwind scheme with periodic boundary conditions in Figure \ref{fig:loop}.  Table \ref{tab:synthetic_results} compares the results to the consensus baseline model of using the singular value decomposition of $L_{cons}$ in Eqn. \eqref{eqn:directed_graph_gp}. 
         \begin{table}[h!]
         \small
\centering
\begin{tabular}{ |c|c|c|c|c|c|c|c| }
\hline
Model & Graph type & $n = 280$ & $n = 325$ & $n = 400$ & $\nu$ & $\kappa$ & $\sigma$\\ \hline
Advection & \multirow{2}{*}{Upwind} & 0.52 & 0.45 & \textbf{0.0005} &0.65 & 8.09 &  7.75  \\
Consensus  &  & \textbf{0.51} & \textbf{0.44} & \textbf{0.0005} &0.65 & 8.29 &  7.77  \\  
\hline

Advection & \multirow{2}{*}{Central}                    & 1.31 & 0.85 & 8.41e-05 & 0.67 & 9.00 & 8.03 \\ 
                 Consensus                    &                     & \textbf{0.97} & \textbf{0.8} & \textbf{8.02e-05} & 0.67 & 9.45 & 8.11 \\
                 \hline
Advection & \multirow{2}{*}{Intersection} & 0.96 & \textbf{0.45} & \textbf{0.0005} & 0.65 & 8.19 & 7.75  \\
Consensus  &  & \textbf{0.52} & 0.46 & \textbf{0.0005} & 0.64 & 8.28 & 7.77  \\ 
\hline
Advection & \multirow{2}{*}{Loop} & \textbf{0.47} & \textbf{0.41} & \textbf{0.00045} & 0.65 & 8.49 & 7.76  \\
Consensus  &  & \textbf{0.47} & \textbf{0.41} & \textbf{0.00045} & 0.65 & 8.49 & 7.76   \\ 
\hline
\end{tabular}
\vspace{.25cm}
 \caption{Comparison of $l_2$ test error on synthetic directed graphs with $n$ nodes and the learned hyperparameters.}
 
\label{tab:synthetic_results}
\end{table}

\vspace{-.25cm}

\begin{figure}[H]
        \centering
     \includegraphics[width=\linewidth]{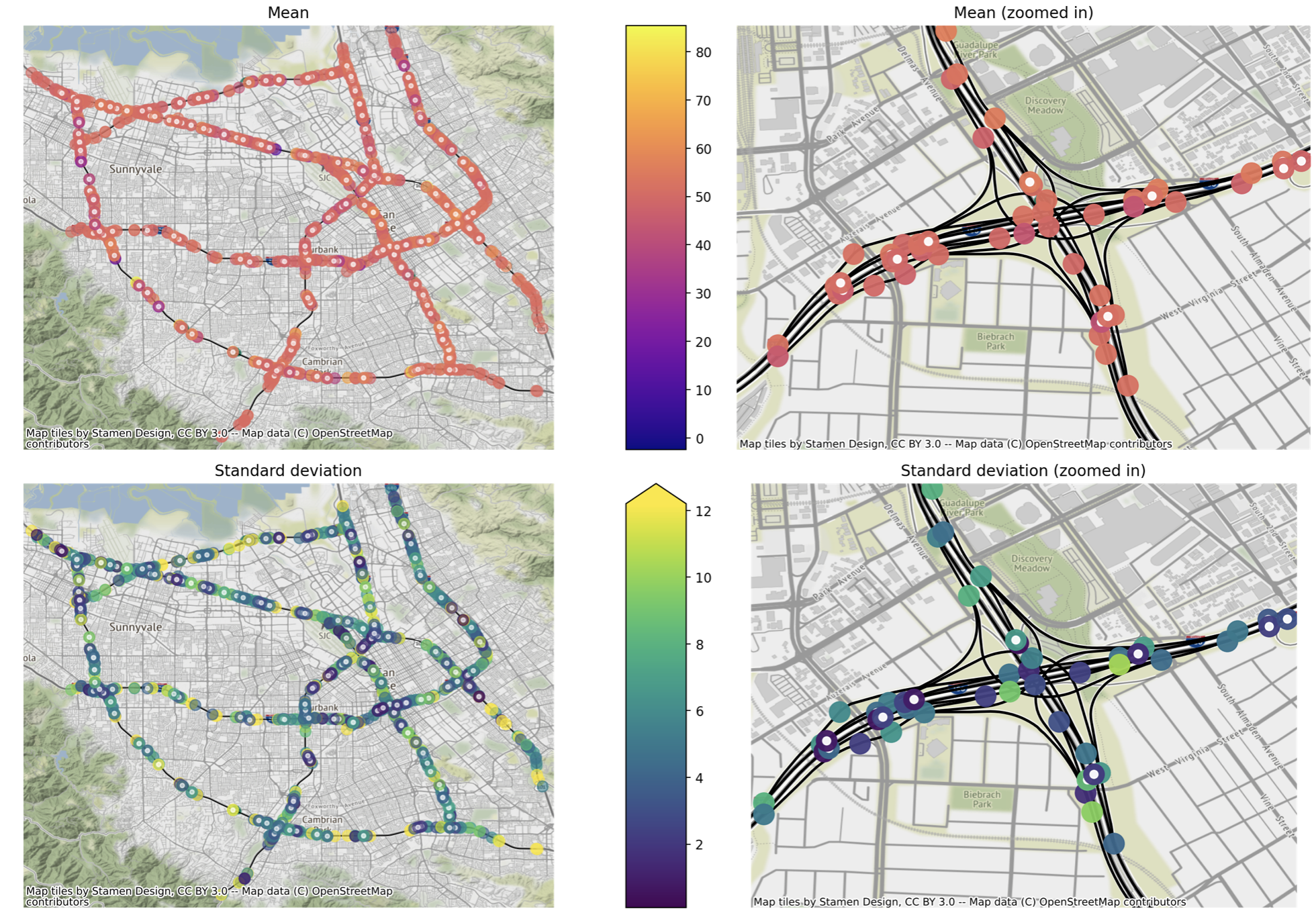}
              \caption{Traffic speed interpolation over a graph of San Jose highways using our DGAMGP method with $\nu = 0.35, \kappa = 1002.8$, $\sigma = 1.14$ and plotting tools from \cite{borovitskiy2021matern}.}
              \label{fig:results-real-data}
      \end{figure}

 \paragraph{Regression results on a real-world traffic graph.}
      We test on the real-world traffic
        data from the California Performance Measurement System \cite{chen2001} with the road network graph from the San Jose highways from Open Street Map \cite{osm2017} at a fixed time. Since our method supports directed graphs, we do not need to convert the raw directed traffic data to an undirected graph as in \cite{borovitskiy2021matern}.  We use the same experimental setup from \cite{borovitskiy2021matern} to generate the train and test data.  
        Figure \ref{fig:results-real-data} shows the resulting predictive mean and standard deviation of the speed on the San Jose highways using the visualization tools from \cite{borovitskiy2021matern}.  
        We notice that the predictive standard deviation along the nodes is relatively small, and is larger on the points that are farther from the sensors.

\section{Conclusions}
In this paper, we propose a novel method DGAMGP to model an advective process on a directed graph and its uncertainties. We show connections between finite differences schemes used to solve the linear advection equation and the graph advection operator $L_{adv}$ employed in our model. 
We explore a regression problem on various graphs, and show that our proposed DGAMGP model performs similarly to other state-of-the-art models.
Future work includes adding a time-varying component to our model, comparing our method to classical numerical methods for solving PDEs, and incorporating the behavior of the non-linear advection equation for traffic modeling.


\bibliography{sample}

\begin{thebibliography}{}

\bibitem[Bakka et~al., 2020]{bakka2020diffusionbased}
Bakka, H., Krainski, E., Bolin, D., Rue, H., and Lindgren, F. (2020).
\newblock The diffusion-based extension of the mat\'ern field to space-time.
\newblock {\em arXiv:2006.04917}.

\bibitem[Borovitskiy et~al., 2021]{borovitskiy2021matern}
Borovitskiy, V., Azangulov, I., Terenin, A., Mostowsky, P., Deisenroth, M., and
  Durrande, N. (2021).
\newblock Mat{é}rn gaussian processes on graphs.
\newblock In Banerjee, A. and Fukumizu, K., editors, {\em Proceedings of The
  24th International Conference on Artificial Intelligence and Statistics},
  volume 130 of {\em Proceedings of Machine Learning Research}, pages
  2593--2601. PMLR.

\bibitem[Chamberlain et~al., 2021]{chamberlain2021grand}
Chamberlain, B., Rowbottom, J., Gorinova, M.~I., Bronstein, M., Webb, S., and
  Rossi, E. (2021).
\newblock Grand: Graph neural diffusion.
\newblock In Meila, M. and Zhang, T., editors, {\em Proceedings of the 38th
  International Conference on Machine Learning}, volume 139 of {\em Proceedings
  of Machine Learning Research}, pages 1407--1418. PMLR.

\bibitem[Chapman and Mesbahi, 2011]{chapman2011}
Chapman, A. and Mesbahi, M. (2011).
\newblock Advection on graphs.
\newblock {\em IEEE Conference on Decision and Control and European Control
  Confereence (CDC-ECC)}, 50:1461--1466.

\bibitem[Chen et~al., 2001]{chen2001}
Chen, C., Petty, K., Skabardonis, A., Varaiya, P., and Jia, Z. (2001).
\newblock Freeway performance measurement system: mining loop detector data.
\newblock {\em Transportation Research Record}, 1748(1):96--102.

\bibitem[Chen et~al., 2018]{driving}
Chen, Y., Sohani, N., and Peng, H. (2018).
\newblock Modelling of uncertain reactive human driving behavior: a
  classification approach.
\newblock In {\em 2018 IEEE Conference on Decision and Control (CDC)}, pages
  3615--3621.

\bibitem[Gardner et~al., 2018]{gardner2018}
Gardner, J., Pleiss, G., Bindel, D., Weinberger, K., and Wilson, A. (2018).
\newblock G{P}ytorch: Blackbox matrix-matrix gaussian process inference with
  gpu acceleration.
\newblock {\em 32nd Conference on Neural Information Processing Systems (NIPS
  2018) arXiv:1809.11165v2}.

\bibitem[Gulian et~al., 2019]{osti_1642956}
Gulian, M., Raissi, M., Perdikaris, P., and Karniadakis, G. (2019).
\newblock Machine learning of space-fractional differential equations, {SIAM}
  {J}ournal on {S}cientific {C}omputing, {V}ol. 41, {N}o. 4, {S}ociety for
  {I}ndustrial and {A}pplied {M}athematics.
\newblock pages A2485--A2509.

\bibitem[Hošek and Volek, 2019]{hosek2019}
Hošek, R. and Volek, J. (2019).
\newblock Discrete advection–diffusion equations on graphs: Maximum principle
  and finite volumes.
\newblock {\em Applied Mathematics and Computation}, 361(C):630--644.

\bibitem[K\"{u}per and Waldherr, 2020]{kuper2020}
K\"{u}per, A. and Waldherr, S. (2020).
\newblock Numerical gaussian process kalman filtering.
\newblock {\em 21st IFAC World Congress}.

\bibitem[LeVeque, 2007]{leveque}
LeVeque, R.~J. (2007).
\newblock {\em Finite Difference Methods for Ordinary and Partial Differential
  Equations: Steady-State and Time-Dependent Problems}.
\newblock SIAM.

\bibitem[Li et~al., 2018]{li2018}
Li, Y., Yu, R., Shahabi, C., and Liu, Y. (2018).
\newblock Diffusion convolutional recurrent neural network: Data-driven traffic
  forecasting.
\newblock {\em International Conference on Learning Representations (ICLR)}.

\bibitem[Liberty et~al., 2020]{sagemaker}
Liberty, E., Karnin, Z., Xiang, B., Rouesnel, L., Coskun, B., Nallapati, R.,
  Delgado, J., Sadoughi, A., Astashonok, Y., Das, P., Balioglu, C.,
  Chakravarty, S., Jha, M., Gautier, P., Arpin, D., Januschowski, T., Flunkert,
  V., Wang, Y., Gasthaus, J., Stella, L., Rangapuram, S., Salinas, D.,
  Schelter, S., and Smola, A. (2020).
\newblock Elastic machine learning algorithms in amazon sagemaker.
\newblock In {\em 2020 ACM SIGMOD International Conference on Management of
  Data, SIGMOD ’20, New York, NY, USA. Association for Computing Machinery.},
  pages 731--737.

\bibitem[Liesen and Parlett, 2008]{liesen2008}
Liesen, J. and Parlett, B.~N. (2008).
\newblock On nonsymmetric saddle point matrices that allow conjugate gradient
  iterations.
\newblock {\em Numer. Math.}, 108:605--624.

\bibitem[Lighthill and Whitham, 1955]{Lighthill1955OnKW}
Lighthill, M. and Whitham, G. (1955).
\newblock On kinematic waves ii. a theory of traffic flow on long crowded
  roads.
\newblock {\em Proceedings of the Royal Society of London. Series A.
  Mathematical and Physical Sciences}, 229:317 -- 345.

\bibitem[OpenStreetMap, 2017]{osm2017}
OpenStreetMap (2017).
\newblock \url {https://www.openstreetmap.org}.

\bibitem[Raissi et~al., 2019]{RAISSI2019686}
Raissi, M., Perdikaris, P., and Karniadakis, G. (2019).
\newblock Physics-informed neural networks: A deep learning framework for
  solving forward and inverse problems involving nonlinear partial differential
  equations.
\newblock {\em Journal of Computational Physics}, 378:686--707.

\bibitem[Rak, 2017]{thesis}
Rak, A. (2017).
\newblock Advection on graphs.
\newblock \url {http://nrs.harvard.edu/urn-3:HUL.InstRepos:38779537}.

\bibitem[Rasmussen and Williams, 2006]{rasmussen2006}
Rasmussen, C. and Williams, C. (2006).
\newblock {\em Gaussian Processes for Machine Learning}.
\newblock MIT Press.

\bibitem[Richards, 1956]{Richards1956}
Richards, P. (1956).
\newblock Shock waves on the highway.
\newblock {\em Operation Res.}, pages 42 -- 51.

\bibitem[Sen and Venkaiah, 1988]{sen1988}
Sen, S. and Venkaiah, V.~C. (1988).
\newblock On symmetrizing a matrix.
\newblock {\em Indian J. pure appl. Math.}, 19(6):554--561.

\bibitem[Solomon, 2015]{solomon2015pde}
Solomon, J. (2015).
\newblock {PDE} approaches to graph analysis.
\newblock {\em ArXiv}, abs/1505.00185.

\end{thebibliography}
\bibliographystyle{apalike}

\appendix

\newpage

\newpage
\section{Upwinding discretizations of linear advection}
  We discretize the 1D linear advection equation with velocity $v$:
    \begin{equation*}
        u_t + vu_x = 0,
    \end{equation*}
using the standard first order upwinding scheme on a simple uniform Cartesian mesh with spatial step size $\Delta x$.  
Then the classical finite difference first-order upwind scheme depends on the sign of $v$.  For flow moving from left to right, $v > 0$, and we have the following semi-discrete discretization \cite{leveque}:
    \begin{equation}
    \begin{aligned}
        & \frac{du_i}{dt} + v \frac{u_i - u_{i-1}}{\Delta x} = 0, \hspace{0.25cm} \text{if} \hspace{0.25cm} v > 0, \\
         & \frac{du_i}{dt} + v \frac{u_{i+1} - u_{i}}{\Delta x} = 0, \hspace{0.25cm} \text{if} \hspace{0.25cm} v < 0. \\
    \end{aligned}
    \label{eqn:first_order_upwind}
    \end{equation}
    Upwinding schemes are useful in the advection case since information is moving from left to right.  The Courant-Friedrichs-Lewy (CFL) condition for stability of the first order upwinding scheme with Forward Euler time-stepping discretization with time step $\Delta t$ is given by:
$$ \Bigl \lvert \frac{v \Delta t}{ \Delta x}\Bigr \rvert \le 1 \iff \Delta t \le \Bigl \lvert \frac{v}{ \Delta x}\Bigr \rvert.$$
    \label{appendix:upwind}

A less diffusive second order upwind scheme is also known as linear upwind differencing (LUD), and is given by:
\begin{equation}
\frac{du_i}{dt} = v \frac{-u_{i-2} + 4u_{i-1}-3u_i}{2 \Delta x}.
\label{eqn:LUD}
\end{equation}
We can show that the scheme is second-order accurate using Taylor expansions. It is designed to be less diffusive because the $u_{xx}$ term from the first-order upwinding scheme cancels.
We have 

\begin{align*}
    \frac{u_{i-2} - 4u_{i-1}+3u_i}{2 \Delta x} &= \frac{1}{2\Delta x} \Bigg[ \Big( u-2 \Delta x u_x + \frac{4\Delta x^2}{2} u_{xx} - \frac {8\Delta x^3}{6} u_{xxx} + \mathcal{O}(\Delta x ^4) \Big) \\
    & + \Big( - 4(u-\Delta x u_x + \frac{\Delta x^2}{2} u_{xx} - \frac {\Delta x^3}{6} u_{xxx} + \mathcal{O}(\Delta x ^4)) \Big)
 + 3u \Bigg] \\
    &= u_x -  \frac {\Delta x^2}{3} u_{xxx} + \mathcal{O}(\Delta x ^4).
\end{align*}
Hence, the scheme is second order accurate with a dispersive $u_{xxx}$ leading error term.

\section{Examples of $L_{adv}$ on balanced graphs resulting in finite difference discretizations of linear advection}
\label{app:finite_diff_graphs}

In addition to the finite difference schemes provided in Section \ref{subsec:finite_diff}, we also provide an example of a non-uniform mesh discretization:
$$ \frac{du_i}{dx} \approx \frac{\frac{4}{3} u_{i+1/2} - u_i - \frac{1}{3} u_{i-1}}{\Delta x}, $$
which results in the following graph, where the in-going and out-going edges from $u_i$:
\begin{center}

\begin{tikzpicture}[->,>=stealth',shorten >=1pt,auto,node distance=2.5cm,
                    thick,main node/.style={circle,draw,font=\sffamily\bfseries}]
                    \centering

  \node[main node] (1) {$u_{i-1}$};
  \node[main node] (2)[right of=1] {$u_{i-1/2}$};
  \node[main node] (3) [right of=2] {$u_{i}$};
  \node[main node] (4) [right of=3] {$u_{i+1/2}$};
  \node[main node] (5) [right of=4] {$u_{i+1}$};

  \path[every node/.style={font=\sffamily\small}]
    (1) edge [bend right] node[below] {$v/3\Delta{x}$} (3)
    (2) edge [bend right] node[below] {$v/3\Delta{x}$} (4)
    (2) edge node [above] {$-4v/3\Delta{x}$} (1) 
    (4) edge node [above] {$-4v/3\Delta{x}$} (3) 
    (3) edge node [above] {$-4v/3\Delta{x}$} (2) 
    (5) edge node [above] {$-4v/3\Delta{x}$} (4) 
    (3) edge [bend right] node[below] {$v/3\Delta{x}$} (5);
\end{tikzpicture}
\end{center}
\vspace{-.25cm}

We can obtain the less diffusive second order upwind scheme (LUD) in \eqref{eqn:LUD} using the following graph: 
\begin{center}
\begin{tikzpicture}[->,>=stealth',shorten >=1pt,auto,node distance=2.5cm,
                    thick,main node/.style={circle,draw,font=\sffamily\bfseries}]
                    \centering

  \node[main node] (1) {$u_{i-2}$};
  \node[main node] (2) [right of=1] {$u_{i-1}$};
  \node[main node] (3) [right of=2] {$u_{i}$};
\node[main node] (4) [right of=3] {$u_{i+1}$};
  \node[main node] (5) [right of=4] {$u_{i+2}$};

  \path[every node/.style={font=\sffamily\small}]
  (1) edge [bend right] node[below] {$-v/2\Delta{x}$} (3)
    (1) edge node[above] {$2v/\Delta{x}$} (2)
    (2) edge node [above] {$2v/\Delta{x}$} (3) 
     (2) edge [bend right] node[below] {$-v/2\Delta{x}$} (4)
    (3) edge node[above] {$2v/\Delta{x}$} (4)
     (3) edge [bend right] node[below] {$-v/2\Delta{x}$} (5)
    (4) edge node [above] {$2v/\Delta{x}$} (5) ;
\end{tikzpicture}
\end{center}



\section{Additional Experiments}
\label{experiments}

       \subsection{Gaussian Process prior results with DGAMGP}
       \label{prior}
       
\begin{figure}[h!]
\centering     
\subfigure[complete graph prior.]{\label{subfig:complete}\includegraphics[width=0.4\textwidth]{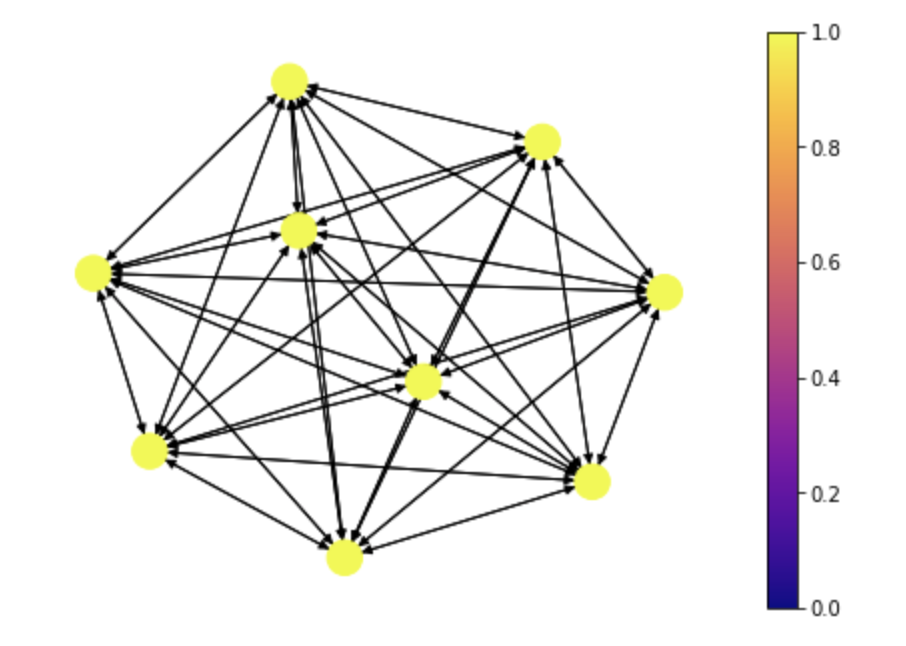}}
\subfigure[star graph prior.]{\label{subfig:star}\includegraphics[width=0.4\textwidth]{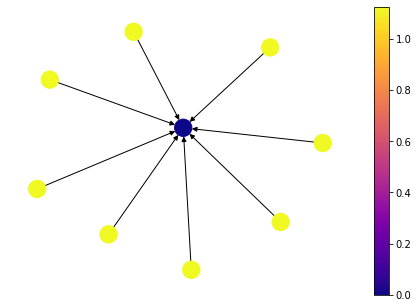}}
\caption{Prior results using DGAMGP obtained using various graphs, and plotting tools from \cite{borovitskiy2021matern}.} 
\label{Velocities}
\end{figure}

        A main property of the Mat\'ern Gaussian Process kernel is that it varies along Riemannian manifolds.   The variance of the kernel is a function of degree, and depends on a complex manner on the graph.
         We show the results generated with a star graph directed towards the center node and a directed complete graph. 
         Figure \ref{subfig:complete} shows that as expected for the complete graph, the nodes have the same variability, since for a random walk starting from any node, there is equal probability to get to another node. For the star graph in Figure \ref{subfig:star}, we observe that the center node has a variability of approximately 0 as starting from any node on the graph, the random walk always ends at the center.
        
\subsection{Convergence Studies}
\begin{figure}[H]
\centering     
\subfigure[Solution of the linear advection equation using \eqref{eqn:graph_adv} ]{\label{fig:upwind}\includegraphics[width=0.49\textwidth]{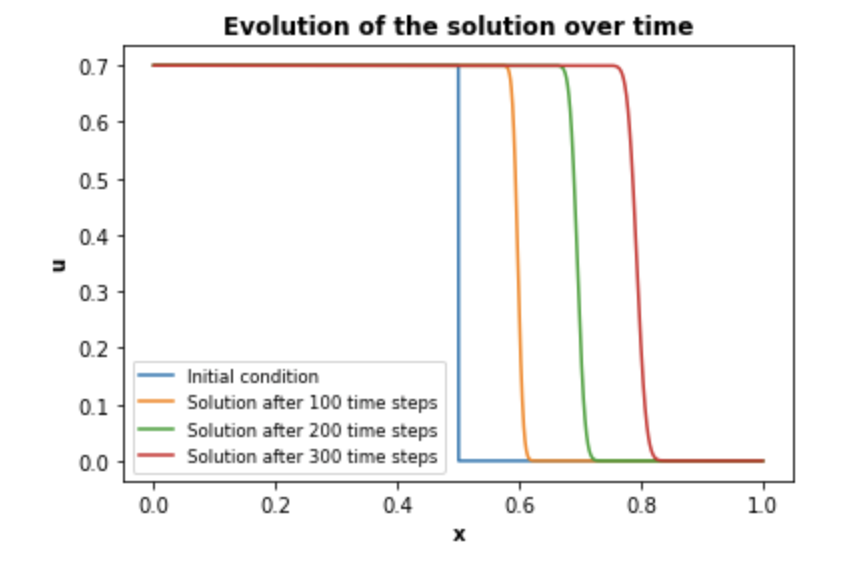}}
\subfigure[Convergence study in a log-log plot]{\label{Conv_anal}\includegraphics[width=0.49\textwidth]{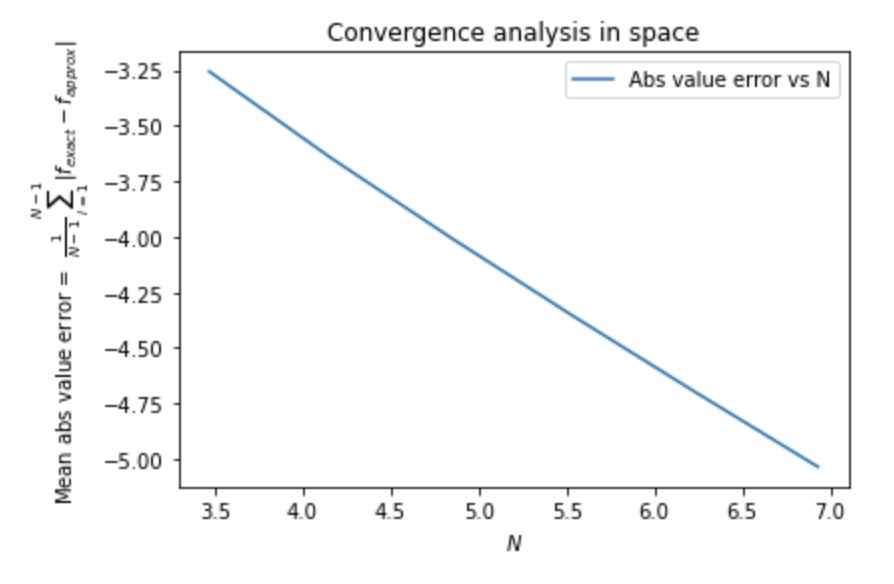}}
\caption{Upwinding solution with RK5 to the linear advection equation over time and corresponding convergence study.}
\label{solve_upwind_conv}
\end{figure}

We conduct a convergence study of applying  $L_{adv}$ on the upwind graph in Figure \ref{fig:graph_upwind}, and show that it has first order convergence matching the performance of the equivalent first order upwind scheme. We use the same initial condition as in Figure \ref{fig:upwind_exact}. We then solve the resulting system of ODEs using the RK5 ODE solver. Figure \ref{fig:upwind} shows the solution at different time steps, and we see how the solution is propagating to the right.  Figure \ref{Conv_anal} shows a loglog plot, where the error is decreasing linearly with a slope of 1 as the number of nodes $n$ is increasing, as expected. 

\end{document}